\documentclass[10pt,conference,a4paper]{IEEEtran}

\usepackage[dvips]{graphicx}
\usepackage{float}
\usepackage{times}
\usepackage{amsmath}
\usepackage{amsfonts}
\renewcommand{\baselinestretch}{1.1}
\IEEEoverridecommandlockouts
\begin{document}
%
\title{A Zero-attracting Quaternion-valued Least Mean Square Algorithm for Sparse System Identification}

\author{\IEEEauthorblockN{Mengdi Jiang, Wei Liu}
\IEEEauthorblockA{Communications Research Group\\ Department of
Electronic and Electrical Engineering\\ University of Sheffield,
UK\\\{mjiang3, w.liu\}@sheffield.ac.uk} \and \IEEEauthorblockN{Yi Li}
\IEEEauthorblockA{School of Mathematics and Statistics\\ University of Sheffield, S3 7RH
UK\\yili@sheffield.ac.uk}\thanks{This work is partially funded by National Grid, UK and will appear in the Proc. of the 9th International Symposium on Communication Systems, Networks and Digital Signal Processing (CSNDSP), Manchester, UK, July 2014 (submitted in March 2014 and accepted on 18 April 2014).}}


\maketitle

\begin{abstract}
Recently, quaternion-valued signal processing has received more and more attention. In this paper, the quaternion-valued sparse system identification problem is studied for the first time and a zero-attracting quaternion-valued least mean square (LMS) algorithm is derived by considering the $l_1$ norm of the quaternion-valued adaptive weight vector. By incorporating the sparsity information of the system into the update process, a faster convergence speed is achieved, as verified by simulation results.
\end{abstract}


\IEEEpeerreviewmaketitle

\section{Introduction}\label{sec:introduction}

In adaptive filtering~\cite{haykin96a}, there is a class of algorithms specifically designed for sparse system identification, where the unknown system only has a few large  coefficients while the remaining ones have a very small amplitude so that they can be ignored without significant effect on the overall performance of the system. A good example of them is the zero-attracting least mean square (ZA-LMS) algorithm proposed in~\cite{chen09b}. This algorithm can achieve a higher convergence speed, and meanwhile, reduce the steady state excess mean square error (MSE). Compared to the classic LMS algorithm~\cite{widrow75a}, the ZA-LMS algorithm introduces an $l_1$ norm in its cost function, which modifies the weight vector update equation with a zero attractor term.

Recently, the hypercomplex concepts have been introduced to solve problems related to three or four-dimensional
signals~\cite{BihanN2004}, such as vector-sensor array signal processing~\cite{zhang13a,liu14e,liu14k}, color image processing~\cite{pei99a} and wind profile prediction~\cite{took09a,liu13j}. As quaternion-valued algorithms can be regarded as an extension of the
complex-valued ones, the adaptive filtering algorithms in complex domain could be extended to the quaternion
domain as well, such as the quaternion-valued LMS (QLMS) algorithm in~\cite{liu14f}.


In this paper, we propose a novel quaternion-valued adaptive algorithm with a sparsity constraint, which is called zero-attracting QLMS (ZA-QLMS) algorithm. The additional constraint is formulated based on the $l_1$ norm. Both the QLMS and ZA-QLMS algorithms can identify an unknown sparse system effectively. However,  a better performance in terms of convergence speed  is achieved by the latter one.

This paper is organized as follows. A review of basic operations in the quaternion domain is provided in Section \ref{sec:quaternion_basic} to facilitate the following derivation of the ZA-QLMS algorithm. The proposed ZA-QLMS algorithm is derived  in Section \ref{sec:ZA-QLMS}. Simulation results are given in Section \ref{sec:simulations}, and conclusions are drawn in Section \ref{sec:conclusions}.

\section{Quaternion-valued Adaptive Filtering}\label{sec:quaternion_basic}

\subsection{Basics of Quaternion}
Quaternion is a non-commutative extension of the complex number, and normally a quaternion consists of one real part and three imaginary parts, denoted by subscripts $a$, $b$, $c$ and $d$, respectively.

For a quaternion number $q$, it can be described as
\begin{equation}
q=q_{a} + (q_{b}i + q_{c}j + q_{d}k),
\label{eq:quaternion}
\end{equation}
where $q_{a}$, $q_{b}$, $q_{c}$, and $q_{d}$ are real-valued~\cite{hamilton1866a,kantor1989a}. For a quaternion, when its real part is zero, it becomes a pure quaternion. In this paper, we consider the conjugate operator of $q$ as $q^{*}=q_a-q_{b}i-q_{c}j-q_{d}k$. The three imaginary units $i$, $j$,  and $k$ satisfy
\begin{eqnarray}
ij=k,~jk=i,~ki=j,\nonumber\\
ijk=i^{2} = j^{2} = k^{2} = -1.
\label{eq:property}
\end{eqnarray}

As a quaternion has the noncommutativity property, in multiplication, the exchange of any two elements in their order will give a different result. For example, we have $ji=-ij$ rather than $ji=ij$.

\subsection{Differentiation with Respect to a Quaternion-valued Vector}\label{sec:Differentiation to a quaternion-valued vector}

To derive the quaternion-valued adaptive algorithm, the starting point is
the general operation of differentiation with respect to a quaternion-valued vector.

At first, we need to give
the definition of differentiation with respect to a quaternion $q$ and its conjugate $q^{*}$. Assume that $f(q)$ is a function of the quaternion variable $q$,
which is expressed as
\begin{equation}
f(q)=f_{a} + f_{b}i + f_{c}j + f_{d}k
\end{equation}
where $f(q)$ is in general quaternion-valued. The definition of $\dfrac{d f(q)}{d q}$ can be expressed as~\cite{mandic2011a,liu14f}
\begin{equation}
\dfrac{d f(q)}{d q}=\frac{1}{4}(\displaystyle\frac{\partial{f(q)}}{\partial q_a}-\displaystyle\frac{\partial{f(q)}}{\partial q_b} i-\displaystyle\frac{\partial{f(q)}}{\partial q_c} j- \displaystyle\frac{\partial{f(q)}}{\partial q_d} k)\;.
\label{eq:general_definition}
\end{equation}

The derivative of $f(q)$ with respect to $q^{*}$ can be defined in a similar way
\begin{equation}
\dfrac{d f(q)}{d q^{*}}=\frac{1}{4}(\displaystyle\frac{\partial{f(q)}}{\partial q_a}+\displaystyle\frac{\partial{f(q)}}{\partial q_b} i+\displaystyle\frac{\partial{f(q)}}{\partial q_c} j+ \displaystyle\frac{\partial{f(q)}}{\partial q_d} k)\;.
\label{eq:conj_general_definition}
\end{equation}

With this definition, we can easily obtain
\begin{equation}
 \frac{\partial q}{\partial q}=1,~\frac{\partial q}{\partial q^{*}}=-\frac{1}{2}\;.
 \end{equation}

Some product rules can be obtained from above formulations, such as the differentiation of quaternion-valued functions to real variables.

Suppose $f(q)$ and $g(q)$ are two quaternion-valued functions of the quaternion variable $q$, and $q_a$ is the real variable. Then we can have the following result
\begin{eqnarray}
\frac{\partial f(q)g(q)}{\partial q_a}&=&\frac{\partial }{\partial q_a}(f_a+if_b+jf_c+kf_d)g\nonumber\\
&=&\frac{\partial f_a g}{\partial q_a}+i\frac{\partial f_b g}{\partial q_a}+j\frac{\partial f_c g}{\partial q_a}+k\frac{\partial f_d g}{\partial q_a}\nonumber\\
&=&(f_a\frac{\partial g}{\partial q_a}+\frac{\partial f_a}{\partial q_a}g)+i(f_b\frac{\partial g}{\partial q_a}+\frac{\partial f_b}{\partial q_a}g)\nonumber\\
&~&+j(f_c\frac{\partial g}{\partial q_a}+\frac{\partial f_c}{\partial q_a}g)+k(f_d\frac{\partial g}{\partial q_a}+\frac{\partial f_d}{\partial q_a}g)\nonumber\\
&=&(f_a+if_b+jf_c+kf_d)\frac{\partial g}{\partial q_a}\nonumber\\
&~&+(\frac{\partial f_a}{\partial q_a}+i\frac{\partial f_b}{\partial q_a}+j\frac{\partial f_c}{\partial q_a}+k\frac{\partial f_d}{\partial q_a})g\nonumber\\
&=&f(q)\frac{\partial g(q)}{\partial q_a}+\frac{\partial f(q)}{\partial q_a} g(q)
\end{eqnarray}

When the quaternion variable $q$ is replaced by a quaternion-valued vector $\textbf{w}$, given by
\begin{equation}
\textbf{w} = [w_1~w_2~\cdots~w_{M}]^{T}
\end{equation}
where $w_m = a_m+b_mi+c_mj+d_mk$, $m=1,...,M$, the differentiation of the function $f(\textbf{w})$ with respect to the vector $\textbf{w}$ can be derived using a combination of (\ref{eq:general_definition}) straightforwardly in the following
\begin{eqnarray}
\dfrac{\partial f}{\partial \textbf{w}}=\frac{1}{4}\left[\begin{matrix}
              \dfrac{\partial f}{\partial a_0}-\dfrac{\partial f}{\partial b_0} i-\dfrac{\partial f}{\partial c_0} j-\dfrac{\partial f}{\partial d_0} k\\
              \dfrac{\partial f}{\partial a_1}-i\dfrac{\partial f}{\partial b_1} i-\dfrac{\partial f}{\partial c_1} j-\dfrac{\partial f}{\partial d_1} k\\
              \vdots \\
              \dfrac{\partial f}{\partial a_{M-1}}-\dfrac{\partial f}{\partial b_{M-1}} i-\dfrac{\partial f}{\partial c_{M-1}} j-\dfrac{\partial f}{\partial d_{M-1}} k
             \end{matrix}\right]
\label{eq:vector_definition}
\end{eqnarray}
Similarly, we define $\dfrac{\partial f}{\partial \textbf{w}^{*}}$ as
\begin{eqnarray}
\dfrac{\partial f}{\partial \textbf{w}^{*}}=\frac{1}{4}\left[\begin{matrix}
              \dfrac{\partial f}{\partial a_0}+\dfrac{\partial f}{\partial b_0} i+\dfrac{\partial f}{\partial c_0} j+\dfrac{\partial f}{\partial d_0} k\\
              \dfrac{\partial f}{\partial a_1}+\dfrac{\partial f}{\partial b_1} i+\dfrac{\partial f}{\partial c_1} j+\dfrac{\partial f}{\partial d_1} k\\
              \vdots \\
              \dfrac{\partial f}{\partial a_{M-1}}+\dfrac{\partial f}{\partial b_{M-1}} i+\dfrac{\partial f}{\partial c_{M-1}} j+ \dfrac{\partial f}{\partial d_{M-1}} k
             \end{matrix}\right]
\label{eq:conj_vector_definition}
\end{eqnarray}
Obviously, when $M=1$, (\ref{eq:vector_definition}) and (\ref{eq:conj_vector_definition}) are reduced to (\ref{eq:general_definition}) and (\ref{eq:conj_general_definition}), respectively.

\section{The Zero-attracting QLMS (ZA-QLMS) Algorithm}\label{sec:ZA-QLMS}

To improve the performance of the LMS algorithm for sparse system identification, the ZA-QLMS algorithm is derived in this section. To achieve this, similar to \cite{chen09b}, in the cost function, we add an $l_1$ norm penalty
term for the quaternion-valued weight vector $\textbf{w}[n]$.

For a standard adaptive filter, the output $y[n]$ and error $e[n]$ can be expressed as
\begin{eqnarray}
y[n]&=&{\textbf{w}^{T}[n]}{\textbf{x}[n]}\\
e[n]&=&d[n]-{\textbf{w}^{T}[n]}{\textbf{x}[n]},
\end{eqnarray}
where $\textbf{w}[n]$ is the adaptive weight vector with a length of $L$, $d[n]$ is the reference signal, $\textbf{x}[n]=[x[n-1], \cdots, x[n-L]]^{T}$ is the input sample vector, and $\{\cdot\}^{T}$ denotes the transpose operation. Moreover, the conjugate form $\textbf{e}^{*}[n]$ of the error signal $e[n]$ is given by
\begin{equation}
e^{*}[n]=d^{*}[n]-{\textbf{x}^{H}[n]}{\textbf{w}^{*}[n]},
\end{equation}

Our proposed cost function with a zero attractor term is given by
\begin{equation}
J_0[n]=e[n]e^{*}[n]+\gamma{\|\textbf{w}[n]\|}_1\;,
\end{equation}
where $\gamma$ is a small constant.

The gradient of the above cost function with respect to $\textbf{w}^{*}[n]$ and $\textbf{w}[n]$ can be respectively expressed as
\begin{equation}
\nabla_{\textbf{w}^{*}}J_0[n]=\frac{\partial {J_0[n]}}{\partial \textbf{w}^{*}}
\label{eq:conj_gradient_cost_function}
\end{equation}
and
\begin{equation}
\nabla_{\textbf{w}}J_0[n]=\frac{\partial {J_0[n]}}{\partial \textbf{w}}
\label{eq:gradient_cost_function}
\end{equation}

From~\cite{mandic2011a,brandwood83a}, we know that
the conjugate gradient gives the maximum steepness direction for the optimization surface.
Therefore, the conjugate gradient $\nabla_{\textbf{w}^{*}}J_0[n]$ will be used to derive the update of the
coefficient weight vector.

Expanding the cost function, we obtain
\begin{eqnarray}
J_0[n]&=&e[n]e^{*}[n]+\gamma{\|\textbf{w}[n]\|}_1\nonumber\\
&=&d[n]d^{*}[n]-d[n]{\textbf{x}^{H}[n]}{\textbf{w}^{*}[n]}-{\textbf{w}^{T}[n]} {\textbf{x}[n]}d^{*}[n]\nonumber\\
&~&+{\textbf{w}^{T}[n]} {\textbf{x}[n]}{\textbf{x}^{H}[n]}{\textbf{w}^{*}[n]}+\gamma{\|\textbf{w}[n]\|}_1\;.
\label{eq:extended_cost_function}
\end{eqnarray}
Furthermore,
\begin{eqnarray}
\frac{\partial {J_0[n]}}{\partial \textbf{w}^{*}}&=&\frac{\partial {(e[n]e^{*}[n]+\gamma{\|\textbf{w}[n]\|}_1)}}{\partial \textbf{w}^{*}}\nonumber\\
&=&\frac{\partial}{\partial \textbf{w}^{*}}(d[n]d^{*}[n]-d[n]{\textbf{x}^{H}[n]}{\textbf{w}^{*}[n]}\nonumber\\
&&-{\textbf{w}^{T}[n]} {\textbf{x}[n]}d^{*}[n]+{\textbf{w}^{T}[n]} {\textbf{x}[n]}{\textbf{x}^{H}[n]}{\textbf{w}^{*}[n]})\nonumber\\
&&+\frac{\partial ({\gamma{\|\textbf{w}[n]\|}_1})}{\partial \textbf{w}^{*}}\;.
\label{eq:extended_gradient_cost_function}
\end{eqnarray}

Details of the derivation process for the gradient are shown in the following
\begin{equation}
\frac {\partial (d[n]d^{*}[n])}{\partial {\textbf{w}^{*}[n]}} = 0
\label{eq:part_1}
\end{equation}
\begin{equation}
\frac {\partial (d[n]{\textbf{x}^{H}[n]}{\textbf{w}^{*}[n]})}{\partial {\textbf{w}^{*}[n]}} = d[n]\textbf{x}^{*}[n]
\label{eq:part_2}
\end{equation}
\begin{equation}
\frac {\partial ({\textbf{w}^{T}[n]}{\textbf{x}[n]}d^{*}[n])}{\partial {\textbf{w}^{*}[n]}} = -\frac{1}{2}d[n]\textbf{x}^{*}[n]
\label{eq:part_3}
\end{equation}
\begin{equation}
\frac{\partial({\textbf{w}^{T}[n]} {\textbf{x}[n]}{\textbf{x}^{H}[n]}{\textbf{w}^{*}[n]})}{\partial {\textbf{w}^{*}[n]}}=\frac{1}{2}{\textbf{w}^{T}[n]}{\textbf{x}[n]}{\textbf{x}^{*}[n]}.
\label{eq:part_4}
\end{equation}
Moreover, the last part of the gradient of cost function is given by
 \begin{equation}
\frac{\partial ({\gamma{\|\textbf{w}[n]\|}_1})}{\partial \textbf{w}^{*}} = \frac{1}{4} \gamma\cdot{sgn(\textbf{w}[n])}\;,
\label{eq:part_5}
\end{equation}
where the symbol $sgn$  is a component-wise sign function that is defined as~\cite{chen09b}
\[
  sgn(x) =
   \begin{cases}
    x/|x| &  x\neq0 \\
    0      &  x = 0
   \end{cases}
 \]

Combining the above results,  the final gradient can be obtained as follows
\begin{equation}
\nabla_{\textbf{w}^{*}}J_0[n]=-\frac{1}{2}e[n]\textbf{x}^{*}[n]+\frac{1}{4}\gamma\cdot{sgn(\textbf{w}[n])}\;.
\end{equation}

With the general update equation for the weight vector
\begin{equation}
\textbf{w}[n+1] = \textbf{w}[n]-\mu \nabla_{\textbf{w}^{*}}J_0[n],
\end{equation}
where $\mu$ is the step size, we arrive at the following update equation for the proposed ZA-QLMS algorithm
\begin{equation}
\textbf{w}[n+1] = \textbf{w}[n]+\mu(e[n]\textbf{x}^{*}[n])-\rho\cdot{sgn(\textbf{w}[n])}\;,
\label{eq:update_weight_vector}
\end{equation}
where  $\rho = \mu \gamma$. The last term represents the zero attractor, which enforces the near-zero coefficients to zero and therefore accelerates the convergence process when majority of the system  coefficients are nearly zero in a sparse system.

Note that equation \eqref{eq:update_weight_vector} will be reduced to the normal QLMS algorithm without the zero attractor term, given by~\cite{liu14f}
\begin{equation}
\textbf{w}[n+1]=\textbf{w}[n]+\mu(e[n]{\textbf{x}}^{*}[n])\;.
\end{equation}

\section{Simulation Results}\label{sec:simulations}
In this part, simulations are performed for sparse system identification using the proposed algorithm in comparison with the QLMS algorithm. Two different sparse systems are considered corresponding to Scenario One and Scenario Two in the following.   The input signal to the adaptive filter is colored and generated by passing a quaternion-valued white gaussian signal through a randomly generated filter. The noise part is  quaternion-valued white Gaussian and added to the output of the unknown sparse system, with a 30dB signal to noise ratio (SNR) for both scenarios.

\subsection{Scenario One}
\begin{figure}[htbp]
\begin{center}
   \includegraphics[width=.49\textwidth]{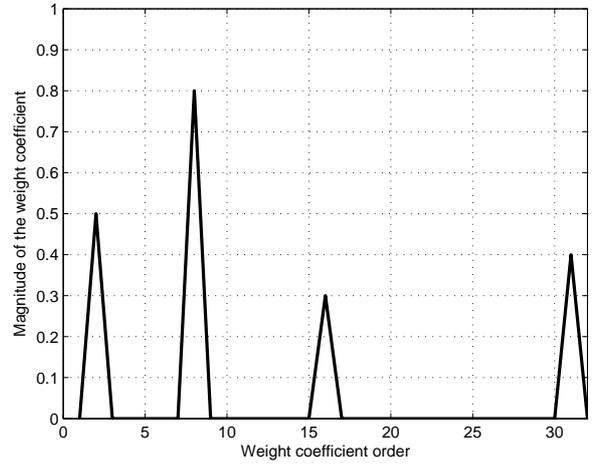}
   \caption{Magnitude of the impulse response of the sparse system.
    \label{fig:w}}
\end{center}
\end{figure}

\begin{figure}[htbp]
\begin{center}
   \includegraphics[width=.49\textwidth]{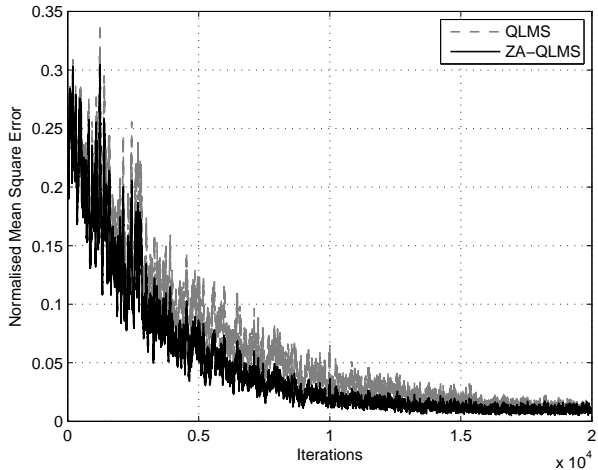}
   \caption{Learning curves for the first scenario.
    \label{fig:learning_curve_1}}
\end{center}
\end{figure}

For the first scenario, the parameters are: the step size $\mu$ is $3\times10^{-7}$; the unknown sparse FIR filter length $L$ is $32$, with $4$ non-zero coefficients at the 2nd, 8th, 16th and 31st taps, and its magnitude of the impulse response is shown in Fig.~\ref{fig:w};  the coefficient of the zero attractor $\rho$ is $5\times10^{-7}$. The learning curve obtained by averaging $100$ runs of the corresponding algorithm is given in Fig.~\ref{fig:learning_curve_1}, where we can see that the ZA-QLMS algorithm has achieved a faster convergence speed than the QLMS algorithm when they both reach a similar steady state.

\subsection{Scenario Two}
\begin{figure}[htbp]
\begin{center}
   \includegraphics[width=.49\textwidth]{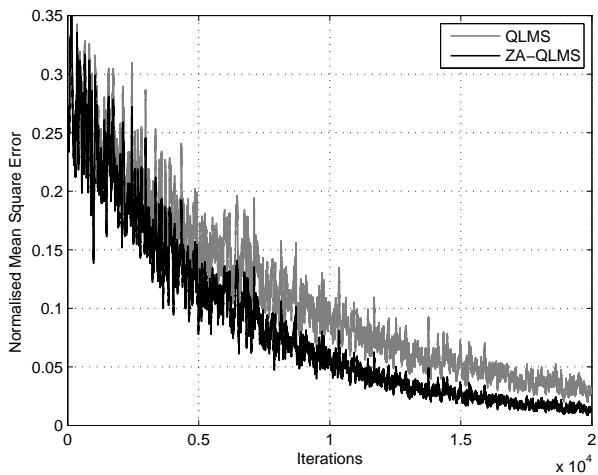}
   \caption{Learning Curves for the second scenario.
    \label{fig:learning_curve_2}}
\end{center}
\end{figure}

For this case, length of the unknown FIR filter is reduced to $16$, still with $4$ active taps. The parameters are: step size $\mu$ is $2\times10^{-7}$ and the value of $\rho$ is $2\times10^{-7}$. The results are shown in  Fig.~\ref{fig:learning_curve_2}. Again we see that the ZA-QLMS algorithm has a faster convergence speed and has even converged to a lower steady state error in this specific scenario. 

\section{Conclusion}\label{sec:conclusions}
In this paper, a quaternion-valued adaptive algorithm has been proposed for more efficient identification of unknown sparse systems. It is derived by introducing an $l_1$ penalty term in the original cost function and the resultant zero-attracting quaternion-valued LMS algorithm can achieve a faster convergence rate by incorporating the sparsity information of the system into the update process. Simulation results have been provided to show the effectiveness of the new algorithm.



\begin{thebibliography}{10}

\bibitem{haykin96a}
S.~Haykin,
\newblock {\em Adaptive Filter Theory},
\newblock Prentice Hall, Englewood Cliffs, New York, 3rd edition, 1996.

\bibitem{chen09b}
Yilun Chen, Yuantao Gu, and Alfred~O Hero,
\newblock ``Sparse \mbox{LMS} for system identification,''
\newblock in {\em Acoustics, Speech and Signal Processing, 2009. ICASSP 2009.
  IEEE International Conference on}. IEEE, 2009, pp. 3125--3128.

\bibitem{widrow75a}
B.~Widrow, J.~McCool, and M.~Ball,
\newblock ``{The Complex LMS Algorithm},''
\newblock {\em Proceedings of the IEEE}, vol. 63, pp. 719--720, August 1975.

\bibitem{BihanN2004}
N.~Le~Bihan and J.~Mars,
\newblock ``Singular value decomposition of quaternion matrices: a new tool for
  vector-sensor signal processing,''
\newblock {\em Signal Processing}, vol. 84, no. 7, pp. 1177--1199, 2004.

\bibitem{zhang13a}
X.R. Zhang, W.~Liu, Y.G. Xu, and Z.W. Liu,
\newblock ``Quaternion-based worst case constrained beamformer based on
  electromagnetic vectoe-sensor arrays,''
\newblock in {\em Proc. IEEE International Conference on Acoustics, Speech, and
  Signal Processing}, Vancouver, Canada, May 2013, pp. 4149--6153.

\bibitem{liu14e}
X.~R. Zhang, W.~Liu, Y.~G. Xu, and Z.~W. Liu,
\newblock ``Quaternion-valued robust adaptive beamformer for electromagnetic
  vector-sensor arrays with worst-case constraint,''
\newblock {\em Signal Processing}, vol. 104, pp. 274--283, November 2014.

\bibitem{liu14k}
M.~B. Hawes and W.~Liu,
\newblock ``A quaternion-valued reweighted minimisation approach to sparse
  vector sensor array design,''
\newblock in {\em Proc. of the International Conference on Digital Signal
  Processing}, Hong Kong, August 2014.

\bibitem{pei99a}
S.C. Pe and C.M. Cheng,
\newblock ``Color image processing by using binary quaternion-moment-preserving
  thresholding technique,''
\newblock {\em Image Processing, IEEE Transactions on}, vol. 8, no. 5, pp.
  614--628, 1999.

\bibitem{took09a}
Clive~Cheong Took and Danilo~P Mandic,
\newblock ``The quaternion \mbox{LMS} algorithm for adaptive filtering of
  hypercomplex processes,''
\newblock {\em IEEE Transactions on Signal Processing}, vol. 57, no. 4, pp.
  1316--1327, 2009.

\bibitem{liu13j}
M.~D. Jiang, W.~Liu, Y.~Li, and X.~R. Zhang,
\newblock ``Frequency-domain quaternion-valued adaptive filtering and its
  application to wind profile prediction,''
\newblock in {\em Proc. of the IEEE TENCON Conference}, Xi'an, China, October
  2013.

\bibitem{liu14f}
M.~D. Jiang, W.~Liu, and Y.~Li,
\newblock ``A general quaternion-valued gradient operator and its applications
  to computational fluid dynamics and adaptive beamforming,''
\newblock in {\em Proc. of the International Conference on Digital Signal
  Processing}, Hong Kong, August 2014.

\bibitem{hamilton1866a}
William~Rowan Hamilton,
\newblock {\em Elements of quaternions},
\newblock Longmans, Green, \& co., 1866.

\bibitem{kantor1989a}
I.~Kantor, A.S. Solodovnikov, and Abe Shenitzer,
\newblock {\em Hypercomplex numbers: an elementary introduction to algebras},
\newblock Springer Verlag, New York, 1989.

\bibitem{mandic2011a}
Danilo~P Mandic, Cyrus Jahanchahi, and Clive~Cheong Took,
\newblock ``A quaternion gradient operator and its applications,''
\newblock {\em IEEE Signal Processing Letters}, vol. 18, no. 1, pp. 47--50,
  2011.

\bibitem{brandwood83a}
DH~Brandwood,
\newblock ``A complex gradient operator and its application in adaptive array
  theory,''
\newblock in {\em IEE Proceedings H (Microwaves, Optics and Antennas)}. IET,
  1983, vol. 130, pp. 11--16.

\end{thebibliography}

\end{document}